\documentclass[11pt]{amsart}   
\topmargin
-1cm 
\newtheorem{thm}{Theorem}[section]
\newtheorem{defn}[thm]{Definition}
\newtheorem{ex}[thm]{Example}
\newtheorem{exs}[thm]{Examples}

\newtheorem{pro}[thm]{Proposition}
\newtheorem{cor}[thm]{Corollary}
\newtheorem{rem}[thm]{Remark}
\newtheorem{rems}[thm]{Remarks}

\begin{document}

\def\sp{\\[0.8ex]}

\noindent

\vspace{0.5in}

\title[Continuous and Pontryagin duality]{Continuous and Pontryagin duality of topological groups}

\author{R. Beattie}
\address{Department of Mathematics and Computer Science, Mount Allison University, Sackville, N.B., Canada, E4L 1E6}
\email{rbeattie@mta.ca}
\author{H.-P. Butzmann}
\address{Fakult\"at f\"ur Mathematik und Informatik, Universit\"at Mannheim, D68131 Mannheim, Germany}
\email{butzmann@math.uni-mannheim.de}

\begin{abstract} 
For Pontryagin's group duality in the setting of locally compact topological Abelian groups, the topology on the character group is the compact open topology. There exist at present two extensions of this theory to topological groups which are not necessarily locally compact. The first, called the Pontryagin dual, retains the compact-open topology. The second, the continuous dual, uses the continuous convergence structure. Both coincide on locally compact topological groups but differ dramatically otherwise. The Pontryagin dual is a topological group while the continuous dual is usually not. On the other hand, the continuous dual is a left adjoint and enjoys many categorical properties which fail for the Pontryagin dual. An examination and comparison of these dualities was initiated in \cite{CMP1}. In this paper we extend this comparison considerably.

\end{abstract}

\maketitle

\section{Preliminaries}
\begin{defn} Let $X$ be a set and, for each $x \in X$, let $\lambda(x)$ be a collection of filters satisfying:
\begin{itemize}
\item[(i)] The ultrafilter $ \dot{x} := \{ A \subseteq X: x \in A \} \in \lambda(x) $,
\item[(ii)] If $\mathcal{F} \in \lambda(x)$ and $\mathcal{G} \in \lambda(x)$, then $\mathcal{F \cap G} \in \lambda(x)$,
\item[(iii)] If $\mathcal{F} \in \lambda(x)$, then $\mathcal{G} \in \lambda(x)$ for all filters $\mathcal{G \supseteq \mathcal{F}}$.
\end{itemize}
\end{defn}

  The totality $\lambda$ of filters $\lambda (x)$ for $x$ in $X$ is called a {\bf convergence structure} for $X$, the pair $(X, \lambda)$ a {\bf convergence space} and filters $\mathcal{F}$ in
$\lambda(x)$ {\bf convergent} to $x$. When there is no ambiguity, $(X, \lambda)$ will be abreviated to $X$. Also, 
we write $ \mathcal{F} \rightarrow x$ instead of $ \mathcal{F} \in
\lambda (x).$ \\A mapping $f: X \rightarrow Y$ between two convergence spaces is \textbf{continuous} if $f(\mathcal{F}) \rightarrow f(x)$ in $Y$ whenever $\mathcal{F} \rightarrow x$ in $X$.

\par Clearly the convergent filters in a topological space satisfy the above axioms so that a topological space is a convergence space. As will become clear, the converse is false. A convergence space $X$ is called \textbf{topological} if there is some topology on the underlying set of $X$ whose convergent filters are precisely those of $X$.

Many properties in topological spaces are easily expressed in terms of filters and so are easily generalized to convergence spaces. For example a convergence space $X$ is called
\begin{itemize}
\item \textbf{Hausdorff} if limits are unique, i.e., if $ \mathcal{F} \rightarrow p$ and $ \mathcal{F} \rightarrow q$ in $X$ implies $ p = q$
\item \textbf{compact} if every ultrafilter on $X$ converges
\item \textbf{locally compact} if it is Hausdorff and each convergent filter contains a compact set.
\end{itemize}

The above are easily seen to be equivalent to the usual definitions for topological spaces.

\begin{defn} Let $G$ be an Abelian group and $\lambda$ a convergence structure on $G$. The pair $(G, \lambda)$ is a \textbf{convergence group} if $\lambda$ is compatible with the group operations, i.e., if the mappings
$$ G \times G \rightarrow G, \quad (x, y) \mapsto x + y$$ and $$ -: G \rightarrow G, \quad x \mapsto -x $$
are continuous. Thus, if $ \mathcal{F} \rightarrow x$ and $ \mathcal{G} \rightarrow y$ in $G$, then the filter $ \mathcal{F} + \mathcal{G}$ generated by $\{A + B: A \in \mathcal{F}, B \in \mathcal{G} \}$ converges to $x + y$ in $G$ and the filter $ -\mathcal{F}$ generated by $\{ -A: A \in \mathcal{F} \}$ converges to $-x$ in $G$.
\end{defn}

\par Every topological group is a convergence group. 

\medskip

We introduce standard notation which will be used. 
\begin{itemize}
\item $\mathbb{T}$ will denote the unit circle in the complex plane.
\item $\mathbb{T}_{+}$ will denote $\{z \in \mathbb{T}: Re(z) \geq 0 \}$. 
\item For a convergence group or topological group $G$, $\Gamma(G)$ will denote the \textbf{character group} of $G$, the group of continuous group homomorphisms from $G$ to $\mathbb{T}$.
\item For $A \subseteq G$, set $A^\circ = \{\varphi \in \Gamma(G): \varphi(x) \in \mathbb{T}_{+}$ for all $x \in A\}$
\item  For $B \subseteq \Gamma(G)$, set
$ B^\diamond = \{x \in G: \varphi (x) \subseteq \mathbb{T}_+$  for
all  $\varphi \in B \} $
\end{itemize}

Throughout, all groups will be assumed to be Abelian. \smallskip

\par Let \textbf{Cgp} and \textbf{Tgp } denote respectively the categories of convergence groups and topological groups, each with continuous group homomorphisms. We single out a particularly important subcategory of \textbf{Cgp} and \textbf{Tgp}, the locally quasi-convex topological groups, introduced in \cite{V} and developed in \cite{Ba}.

\begin{defn} Let $G$ be a topological group.\sp
(i) A set $A \subseteq G$ is called \textbf{quasi-convex } if
$A = A^{\circ\diamond}$.\sp
(ii) G is called \textbf{locally quasi-convex} if it has a zero neighbourhood basis consisting of quasi-convex sets.            
\end{defn}

\begin{thm} (\cite{A1}, 6.8, 6.17) Subgroups, products, projective limits, Hausdorff coproducts,  and Hausdorff completions in \textbf{Tgp} of locally quasi-convex groups are locally quasi-convex.
\end{thm}

If now \textbf{Qgp }denotes the category of locally quasi-convex topological groups, then we have $$ \textbf{Qgp } \subseteq \textbf{Tgp} \subseteq \textbf{Cgp}$$ and each subcategory is reflective. Finally, we denote by \textbf{HQgp} etc. the Hausdorff modifications of the above categories.

\smallskip

Locally quasi-convex topological groups are the group analogue of the locally convex topological vector spaces and duality arguments are most useful in this setting. For continuous duality, the locally quasi-convex groups are precisely those embedded into their bidual (\cite{BB}, 8.4.7); for Pontryagin duality, every dual group is already locally quasi-convex (\cite{A1}, 6.6).

\smallskip

Continuous convergence is usually defined in a very general function space setting (see, e.g., \cite{Bi}, \cite{BB}). Here, for simplicity, we restrict its definition to character groups.

Let $G$ be a convergence or topological group and $\Gamma(G)$ its character group. The {\bf continuous convergence} structure on $\Gamma (G)$ is the coarsest convergence structure on $\Gamma (G)$ making the evaluation mapping
$$\omega: \Gamma(G) \times G \to \mathbb{T} \quad (\varphi,x) \mapsto \varphi(x)$$
continuous. A filter $ \Phi \rightarrow \varphi$ in this convergence structure if, whenever $ \mathcal{F} \rightarrow x$ in $G$, the filter
$ \Phi({\mathcal F}) = \omega( \Phi \times \mathcal{F})$ converges to $\omega(\varphi,x) = \varphi(x)$ in $\mathbb{T}$. This is compatible with the group structure on the character group $\Gamma (G)$ and the resulting convergence group is denoted by $\Gamma_{c} (G)$.

\section{The dual as an adjoint}

When the character group $\Gamma(G)$ of a topological group $G$ carries the continuous convergence structure, it is called the \textbf{continuous dual} (\textbf{c-dual}) of $G$ and denoted by $\Gamma_c(G)$. Likewise, when $\Gamma(G)$ carries the compact-open topology, it is called the compact-open or \textbf{Pontryagin dual} (\textbf{P-dual}) of $G$. Historically this has been denoted $G^{\wedge}$ and we will respect this notation.
\smallskip

The following important result is typical of the behaviour of continuous convergence in various settings (see e.g. \cite{BB}, p. 58).  A proof can be found in (\cite{BB08}, Theorem 3.1), an outline thereof is given in the proof of Theorem 4.1 in \cite{ATC1}. 

\begin{pro} The functor $\Gamma_c : \textbf{Cgp} \rightarrow \textbf{Cgp}^{op} $ is left adjoint to $\Gamma_c^{op}: \textbf{Cgp}^{op} \rightarrow \textbf{Cgp}$. Here $\textbf{Cgp}^{op}$ denotes the opposite category of $\textbf{Cgp}$.
\end{pro}

\par As a left adjoint, $\Gamma_c : \textbf{Cgp} \rightarrow \textbf{Cgp}^{op} $ preserves colimits and so transforms colimits in \textbf{Cgp} to limits in \textbf{Cgp}. As a result, we get:

\begin{cor}\hspace*{0em}\\
(i) If $q: G \rightarrow Q$ is a quotient mapping between convergence groups, then $q^*: \Gamma_c(Q) \rightarrow \Gamma_c(G)$ is an embedding (\cite{Bu1}, Kap. 4, Lemma 9, \cite{BC}, 3.2).\sp
(ii)  Let $(G_i)_{i \in I}$ be a family of convergence groups and $\varepsilon_j: G_j \rightarrow \bigoplus_{i \in I} G_i$ the natural embeddings. Then the mapping $$ \Gamma_c(\bigoplus_{i \in I} G_i) \rightarrow \prod_{i \in I} \Gamma_c(G_i)$$ which maps $\varphi \in \Gamma(\oplus_{i \in I} G_i)$ to $(\varphi \circ \varepsilon_i)_{i \in I}$ is an isomorphism (\cite{Bu3}, 2.2, \cite{BB}, 8.1.19).\sp
(iii) If $(G, (e_i))$ is the inductive limit of an inductive system of convergence groups $G_i$, then $(\Gamma_c G , (e_i^*))$ is isomorphic to the projective limit of the dual system, i.e., $$ \Gamma_c(ind \; G_i) \cong proj \; \Gamma_c(G_i)$$
(\cite{ATC1}, 4.1).
\end{cor}

As well, $\Gamma_c$ enjoys additional duality properties which are not consequences of being a left adjoint:

\begin{pro} \label{extra props}\hspace*{0em}\sp
(i) Let $(G_i)_{i \in I}$ be a family of convergence groups and let $e_j: G_j \rightarrow \prod_{i \in I} G_i$ be the natural embeddings. Then the mapping $$ \Gamma_c(\prod_{i \in I} G_i ) \rightarrow \bigoplus_{i \in I} \Gamma_c(G_i)$$ which maps $\varphi \in \Gamma(\prod_{i \in I} G_i)$ to $(\varphi \circ e_i)$ is an isomorphism  (\cite{Bu3}, 2.3, \cite{BB}, 8.1.18). \sp
(ii)  Let $(G, (p_i))$ be the reduced projective limit of a projective system of topological groups. Then $\Gamma_c(G)$ is isomorphic to the inductive limit of the dual system, i.e., $$ \Gamma_c(proj \; G_i) \cong ind \; \Gamma_c(G_i)$$
(\cite{BB08}, 3.7).
\end{pro}

\begin{rems}\hspace*{0em}\sp
(i) In part (ii) of the above proposition, it is important that the $G_i$ be topological groups. The claim is false for the more general case of convergence groups (\cite{BB93}, Ex. 3).\sp
(ii) In (\cite{ATC1}, 4.4) part (ii) of the above proposition is proved for the special case that all $G_i$'s are nuclear.

(iii) In general, the continuous dual does not carry subgroups to quotients. The dual of the embedding need not even be surjective. However, subgroups of nuclear groups are carried to quotients (\cite{BB}, 8.4.10).
\end{rems}

The Pontryagin dual is not a left adjoint either in \textbf{Tgp} or \textbf{Qgp}. Quotients in \textbf{Qgp} are not carried to embeddings (\cite{Ba}, 17.7) and inductive limits are not carried to projective limits (\cite{BB08}, Example 4.3). Also, as is the case with the continuous dual, the Pontryagin dual does not carry subspaces to quotients. The P-dual actually is a left adjoint when restricted to the subcategory of c-groups (see \ref{c-group} for the definition). In \cite{ATC}, this is proven for the class of $c^\infty-$ groups.

\smallskip

The following result was proven by Kaplan ([25]) for P-reflexive topological groups. The general case is due to (\cite{Ba}, 14.11).

\begin{pro} \label{p-dual} \hspace*{0em}\sp
(i)  If $(G_i)_{i \in I}$ is a family of topological groups, then the mapping $$ (\prod_{i \in I} G_i )^{\wedge} \rightarrow \bigoplus_{i \in I} G_i^{\wedge} $$ which maps $\varphi \in \Gamma(\prod_{i \in I} G_i)$ to $(\varphi \circ e_i)$ is an isomorphism. \sp
(ii) If $(G_i)_{i \in I}$ is a family of locally quasi-convex groups, then the mapping $$ (\bigoplus_{i \in I} G_i)^{\wedge} \rightarrow \prod_{i \in I} G_i^{\wedge}$$ which maps $\varphi \in \Gamma(\oplus_{i \in I} G_i)$ to $(\varphi \circ \varepsilon_i)_{i \in I}$ is an isomorphism.
\end{pro}

Some remarks concerning \ref{p-dual} are necessary. Part (i) holds when the coproduct carries, not the natural coproduct topology in \textbf{Tgp}, but rather the "asterisk" topology (\cite{K}, \cite{N}). It is false otherwise.
However, it was shown in \cite{CD} that the "asterisk" coproduct is the reflection in \textbf{Qgp} of the topological coproduct in \textbf{Tgp} and is precisely the topological coproduct in \textbf{Qgp}.

\smallskip

The situation for the Pontryagin dual of projective limits remains unclear. While a result as general as \ref{extra props} appears unlikely, we have the following:

\begin{pro} (\cite{GG}, Proposition 8.5) If $G$ is the reduced projective limit of Pontryagin reflexive topological groups $(G_i)$, then $G^{\wedge}$ is the inductive limit in \mbox{\textbf{HQgp} of $G_i^{\wedge}$.}
\end{pro}

In (\cite{ATC}, 4.3)  this result is proven for the case of a sequence of metrizable P-reflexive groups.

\section{Embeddedness and Reflexivity}

For a topological group $G$, both the continuous dual $\Gamma_c(G)$ and the Pontryagin dual $G^{\wedge}$ have their own duals and so one has two notions of the bidual, the continuous bidual (\textbf{c-bidual}) $\Gamma_c \Gamma_c(G)$ and the Pontryagin bidual (\textbf{P-bidual}) $G^{\wedge \wedge}$.

The canonical mapping from a topological group $G$ into its bidual has a different notation depending on the duality. The mappings $\kappa_G: G \rightarrow \Gamma_c \Gamma_c(G)$ and $\alpha_G: G \rightarrow G^{\wedge\wedge}$ differ only in the codomain. Both are defined by $x \mapsto \widehat{x}$ where $\widehat{x}(\varphi) = \varphi(x)$ for all $\varphi \in \Gamma(G)$.

\smallskip

There are two stages to the reflexivity of a topological group $G$: the canonical mapping must be surjective and an embedding. 

A convergence or topological group $G$ is called continuously embedded (\textbf{c-embedded}) if $\kappa_G$ is an isomorphism onto its range and continuously reflexive (\textbf{c-reflexive}) if $\kappa_G$ is an isomorphism. In the same way, $G$ is called Pontryagin embedded (\textbf{P-embedded}) if $\alpha_G$ is an isomorphism onto its range and Pontryagin reflexive (\textbf{P-reflexive}) if $\alpha_G$ is an isomorphism.

\begin{rems}\hspace*{0em}\sp
(i)  For a topological group $G$, the evaluation mapping $\omega: \Gamma_c(G) \times G \rightarrow \mathbb{T}$ is continuous. The continuous dual $\Gamma_c(G)$ of a topological group $G$ is a locally compact convergence group. The polars $U^\circ$ of zero neighbourhoods in $G$ are compact in $\Gamma_c(G)$ and, in fact, form a cofinal system of the compact sets. \\In general, $\Gamma_c(G)$ is not topological. It is, however, when $G$ is locally compact and, in this case, the continuous dual and Pontryagin dual coincide. \\
In particular, the c-bidual of a topological group is also a topological group \mbox{(cf. \cite{BCMPT}, 3.5).}\sp
(ii)  For a topological group $G$, the Pontryagin dual $G^{\wedge}$ is a locally quasi-convex group (\cite{A1}, 5.9). It is well known that, in general, the evaluation mapping $\omega: G^{\wedge} \times G \rightarrow \mathbb{T}$ is not continuous. In fact, for P-reflexive groups, this is the case if and only if $G$ is locally compact \cite{MP}.
\end{rems}

A consequence of the above remark is that both the continuous duality and the Pontryagin duality are true extensions of Pontryagin's theory for locally compact groups.

\begin{thm} \label{alpha open} Let $G$ be a topological group.\sp
(i)  The mapping $\kappa_G: G \rightarrow \Gamma_c \Gamma_c(G)$ is continuous (\cite{Bu1}, Kap. 3, Lemma 1). \sp
(ii)  $G$ is c-embedded if and only if $G$ is Hausdorff and locally quasi-convex (\cite{BB}, 8.1.3, 8.4.7, \cite{Bu2}, 4.3).    \sp
(iii) The mapping $\alpha_G: G \rightarrow G^{\wedge\wedge}$ is not continuous in general. If  $G$ is Hausdorff and locally quasi-convex, then $\alpha_G$ is open (\cite{A1}, 6.10).
\end{thm}

Thus, for Hausdorff locally quasi-convex groups, P-embeddedness is equivalent to the continuity of $\alpha_G$. 

\begin{thm} \label{c-emb left adj} The c-embedded groups are a reflective subcategory of \mbox{\textbf{Cgp}, \textbf{Tgp} or \textbf{Qgp}}.
\end{thm}

Proof The mapping which takes $G$ to the initial topology with respect to $\kappa_G$ is easily seen to be left adjoint to the embedding functor.

\begin{thm}  Subgroups, products, coproducts and projective limits in \textbf{Cgp} of c-embedded topological groups are c-embedded.
\end{thm}

Proof Subgroups, products and projective limits follow from the left adjointness in \ref{c-emb left adj}. Coproducts are a straightforward computation.

\begin{exs} \hspace*{0em}\sp
(i) Locally convex topological vector spaces are c-embedded (\cite{Bu4}, Satz 7,  \cite{Bi},  The- \linebreak orem 90). \sp 
(ii) The c-dual $\Gamma_c(G)$ of any convergence group $G$ is c-embedded (\cite{BB}, 8.1.12).
\end{exs}

\begin{thm} \hspace*{0em}\sp
(i) Products and coproducts of c-reflexive convergence groups are c-reflexive (\cite{Bu3}, 2.4, \cite{BB}, 8.1.20). \sp 
(ii) The reduced projective limit of c-reflexive topological groups is c-reflexive \mbox{(\cite{BB}, 8.4.14).}\sp
(iii) Dually closed and dually embedded subgroups of c-reflexive convergence groups are c-reflexive (\cite{BB}, 8.3.8).
\end{thm}

The following are some important examples of c-reflexive groups:

\begin{exs}\hspace*{0em}\sp
(i) For any convergence space $X$, the function space $\mathcal{C}_c(X)$ is c-reflexive (\cite{Bu4}, Satz 4,  \cite{Bi}, Theorem 89).\sp
(ii) For any locally compact convergence space $X$, the group $\mathcal{C}_c(X, \mathbb{T})$ of continuous unimodular functions is c-reflexive (\cite{Bu2}, Theorem 3). \sp
(iii) Every complete locally convex topological vector group, in particular every complete locally convex topological vector space is c-reflexive (\cite{BB}, 8.4.16).\sp
(iv) Every complete nuclear group is c-reflexive (\cite{BB}, 8.4.19).
\end{exs}

The continuity of $\alpha_G$ is such a serious problem that groups that have this property merit a special name:

\begin{defn} \label{c-group}
A topological group is called a \textbf{c-group} if the mapping $\alpha_G: G \rightarrow G^{\wedge\wedge}$ is continuous.
\end{defn}

Note that this is equivalent to the compact sets of $G^{\wedge}$ being equicontinuous \mbox{(\cite{A1}, 5.10).}

\begin{rem} The terminology "c-group" is historic but is most unfortunate if one is dealing simultaneously with both c-duality and P-duality. For, despite its name, c-group has nothing to do with continuous duality and everything to do with Pontryagin duality.
\end{rem}

\begin{ex} \label{metr is c}Each metrizable group is a  c-group (\cite{Ba}, 14.4, \cite{A1}, 5.12). 
\end{ex}

\begin{thm} \label{c-groups coreflective} The c-groups form a coreflective subcategory of \textbf{Tgp}.
\end{thm}

Proof For a fixed $G \in \textbf{Tgp}$ consider the collection of all c-groups on the underlying group of $G$, a non-empty collection. Assign to $G$ the final topological group topology from these c-groups, say $G_0$. It is straightforward to show that $G_0$ is a c-group and that the functor $G \mapsto G_0$ is right adjoint to the embedding functor.

\begin{cor} Products, coproducts, quotients and inductive limits in \textbf{Tgp} of c-groups are c-groups.
\end{cor}

Proof Coproducts, quotients and inductive limits follow from the right adjointness in \ref{c-groups coreflective}. Products are a straightforward computation.

\medskip

It seems to be an open problem whether or not the P-dual $G^{\wedge}$ of a c-group is also a c-group.

The following follows immediately from \ref{alpha open}:

\begin{pro} A Hausdorff locally quasi-convex group is P-embedded if and only if it is a c-group.
\end{pro}

In light of \ref{c-emb left adj}, the following represents one of the most spectacular differences in the two dualities:

\begin{thm} \label{P-emb right adj} The P-embedded topological groups are a coreflective subcategory
\linebreak  of \textbf{HQgp}.
\end{thm}

Proof For a fixed $G \in \textbf{HQgp}$, obtain the c-group $G_0$ exactly as in \ref{c-groups coreflective}. Then the locally quasi-convex modification of $G_0$, say $q(G_0)$, is also a c-group and is Hausdorff since it is finer than $G_0$. Thus $q(G_0)$ is a c-group, locally quasi-convex and Hausdorff, hence P-embedded. The functor $G \mapsto q(G_0)$ is right adjoint to the embedding functor.

\begin{thm} Products, coproducts, quotients and inductive limits in \textbf{Qgp} of P-embedded groups are P-embedded as well.
\end{thm}

Proof Coproducts, quotients and inductive limits follow from the right adjointness in \ref{P-emb right adj}. Products are a straightforward computation.

\begin{exs} \hspace*{0em}\sp
(i)  Every locally compact topological group is P-embedded.\sp
(ii) Every locally quasi-convex metrizable group is P-embedded (\ref{metr is c} and \ref{alpha open}(iii)). \sp
(iii)  Not every nuclear group is P-embedded (\cite{BB}, 8.6.3).\sp
(iv) Not every locally convex topological vector space is P-embedded.
 For example, a Banach space carrying its weak topology is not P-embedded (\cite{A1}, 8.23).
\end{exs}

\begin{thm} \hspace*{0em}\sp
(i) Products and coproducts in \textbf{Qgp} of P-reflexive groups are P-reflexive (\cite{K}, \cite{CD}, \linebreak Theorem 21). \sp
(ii)  The reduced projective limit of a sequence of metrizable P-reflexive groups is a metrizable P-reflexive group (\cite{ATC}, 3.2).
\end{thm}

\begin{exs}\hspace*{0em}\sp
(i)  Every reflexive locally convex topological vector space is P-reflexive \mbox{(\cite{S}, Theorem 1).} \sp
(ii)   Every Banach space is P-reflexive (\cite{S}, Theorem 2). In fact, every Fr\'{e}chet space is P-reflexive.
\end{exs}

\section{Topological groups determined by subgroups}

\par Let $G$ be a topological group and $H$ a subgroup. When do the duals of $G$ and $H$ coincide, i.e., when is $\Gamma_c(G) = \Gamma_c(H)$ or when is $ G^{\wedge } = H^{\wedge}$ ?  In this situation $H$ is almost always assumed to be dense in $G$ although there is no need to do so.

For continuous duality, the situation for dense subgroups of topological groups is particularly simple:

\begin{pro}Let $G$ be a topological group and $\widetilde{G}$ its completion. Then $\Gamma_c(G)$ and $\Gamma_c(\widetilde{G})$ are isomorphic  (\cite{BB}, 8.4.4).
\end{pro}

\begin{cor} Let $G$ be a topological group and $H$ a dense subgroup. Then
$\Gamma_c(G) = \Gamma_c(H)$ (cf. \cite{CMP2}, Theorem 3).
\end{cor}

Proof  $G$ and $H$ have the same (topological group) completion.

\medskip

Slightly more is true:

\begin{cor} \label{dense subgroups} Let $M$ and $N$ be dense subgroups of a topological group $G$. Then $\Gamma_c(M) = \Gamma_c(N)$.
\end{cor}

\par The situation for the Pontryagin duals is much more complicated and has become a much studied problem (see, e.g., \cite{HMTA}, \cite{CRTA}, \cite{C}).

\begin{defn} Let $G$ be a topological group. If, for every dense subgroup $H$ of $G$, $H^{\wedge} = G^{\wedge}$, $G$ is called \textbf{determined}.
\end{defn}

The strongest positive result on determined groups is the following:

\begin{pro}Every metrizable group $G$ is determined  (\cite{C}, Theorem 2, \cite{A1}, 4.5).
\end{pro}

It is astonishing that compact groups need not be determined. In fact:

\begin{pro} A compact group is determined if and only if it is metrizable \linebreak (\cite{HMTA}, 5.11).
\end{pro}

For locally compact groups, the following has only recently been proven (\cite{CMP2}, Proposition 7).

\begin{pro} Let $G$ be a locally compact group and $H$ a dense subgroup. Then $ G^{\wedge } = H^{\wedge}$ if and only if $H$ is P-embedded and $H^{\wedge}$ is locally compact.
\end{pro}

\begin{rem}
 A subgroup $H$ of $G$, need not be dense in $G$ in order to determine $G$. This is true in both dualities. In an example given in \cite{A1}, if $G$ := $L^2[0,1]$ and $H$ := $L_\mathbb{Z}^2[0, 1]$, then $H$ is a complete metrizable group, a closed subgroup of $G$ but $H^{\wedge} = G^{\wedge}$ and $\Gamma_c(H) = \Gamma_c(G)$.
\end{rem}

\section{Distinguishing groups by their duals}

\par The question we examine is: to what extent is a topological group distinguished by its dual? Do different topological groups have different duals or, equivalently, if the duals of two topological groups are equal, must the groups themselves coincide? The two dualities behave very differently in this regard.

\begin{pro} Let $G_1$ and $G_2$ be two c-embedded convergence groups on the same underlying group $G$. If $\Gamma_c(G_1) = \Gamma_c(G_2)$, then $G_1 = G_2$.
\end{pro}

Proof Both $G_1$ and $G_2$ carry the structures of their c-biduals. But these are equal.

\begin{cor} Let $G_1$ and $G_2$ be two Hausdorff locally quasi-con-vex topological groups on the same underlying group $G$. If $\Gamma_c(G_1) = \Gamma_c(G_2)$, then $G_1 = G_2$.
\end{cor}

Proof By \ref{alpha open}, $G_1$ and $G_2$ are c-embedded.

\smallskip

Thus different Hausdorff locally quasi-convex group topologies on a group $G$ have different continuous duals.

This is, of course, not the case with Pontryagin duality. If groups $G_1$ and $G_2$ have the same characters and the same compact sets, then their Pontryagin duals coincide. This is not an unusual phenomenon.

One common instance where this occurs is when a group $G$ respects compactness:

\begin{defn} Let $G$ be a locally quasi-convex group and let $G_b$ be  $G$ carrying its weak (Bohr) topology. $G$ is said to \textbf{respect compactness} (or satisfy Glicksberg's property) if the compact sets of $G$ and $G_b$ coincide.
\end{defn}

\par There exist many important examples of topological groups which respect compactness. For example, locally compact, nuclear and Schwartz groups all have this property (\cite{G}, Theorem 1.2, \cite{BaMP}, Theorem, \cite{A2}, Theorem 4.4).

It is clear that, if $G$ respects compactness, then $ G^{\wedge } = G_b^{\wedge}$. In fact there might well be a whole spectrum of (locally quasi-convex) group topologies on the underlying group of $G$ all having the same Pontryagin dual. In such a case, if one of these topological groups is P-reflexive, it follows that none of the others can be. In fact, if one of these topological groups is P-embedded, none of the others can be.

\section{Relations between the dualities}

Despite the very different behaviour of the two duals, on occasion the continuous dual and Pontryagin dual are very closely related. The following relate the two dualities.

\begin{pro} For any convergence group $G$, the identity mapping $id: \Gamma_c(G) \rightarrow G^{\wedge}$ is continuous. It is a topological isomorphism if $G$ is locally compact (\cite{BB}, 8.1.1, 8.1.2).
\end{pro}

\begin{pro}If $G$ is c-reflexive, then the P-dual $G^{\wedge}$ is the reflection in \textbf{Qgp} of the c-dual $\Gamma_c(G)$  (\cite{BB08}, 4.6).
\end{pro}

\pagebreak

\begin{pro}\hspace*{0em}\sp
(i) For any topological group $G$, the identity mapping $id:  G^{\wedge \wedge} \rightarrow \Gamma_c \Gamma_c(G)$ is a continuous homomorphism between topological groups  (\cite{CMP1}, Theorem 1, \cite{BB}, 8.6.4).\sp
(ii) If $G$ is P-embedded (even if it is a c-group), then $G^{\wedge \wedge}$ is a topological subgroup of $\Gamma_c \Gamma_c(G)$.
\end{pro}

Proof of (i) Because the identity mapping $\Gamma_c(G) \rightarrow G^{\wedge}$ is continuous, it follows that the identity $G^{\wedge\wedge} \rightarrow \Gamma_c(G)^{\wedge}$ is also continuous. But $\Gamma_c(G)^{\wedge} = \Gamma_c \Gamma_c(G)$ since $\Gamma_c(G)$ is locally compact.

\begin{cor} If $G$ is P-embedded (even if it is a c-group) and c-reflexive, it is also P-reflexive.
\end{cor}

\begin{cor} If $G$ is any metrizable group, then $G$ is c-reflexive if and only if it is P-reflexive  (\cite{C}, Theorem 3).
\end{cor}

\end{document}